\documentclass{amsart}
\usepackage{amsmath}
\usepackage{amssymb}
\usepackage{amsfonts}

\theoremstyle{plain}

\newtheorem{algorithm}{Algorithm}

\newtheorem{conclusion}{Conclusion}

\newtheorem{definition}{Definition}
\newtheorem{example}{Example}

\newtheorem{lemma}{Lemma}

\newtheorem{proposition}{Proposition}

\numberwithin{equation}{section}
\input{tcilatex}

\begin{document}
\title{On Bipolar Soft Sets}
\author{Muhammad Shabir}
\address{Department of Mathematics, Quaid-i-Azam University, Islamabad,
Pakistan 46000}
\email{mshabirbhatti@yahoo.co.uk}
\author{Munazza Naz}
\address{Department of Mathematics, Fatima Jinnah Women University, The
Mall, Rawalpindi, Pakistan 46000}
\email{munazzanaz@yahoo.com}
\date{February 18, 2012}

\begin{abstract}
We have studied the concept of bipolarity of information in the soft sets.
We have defined bipolar soft sets and basic operations of union,
intersection and complementation for bipolar soft sets. Examples of bipolar
soft sets and an application of bipolar soft sets in a decision making
problem with general algorithms have also been presented at the end.
\end{abstract}

\maketitle

\section{Introduction}

The idea of soft set and its applications was introduced by Molodtsov $^{%
\text{\cite{[1]}}}$in 1999. The concept was motivated by the limitations on
parameterization process faced by the other theories while handling
uncertain and fuzzy data. Molodtsov presented a new structure namely soft
set accommodating the parameterization problem. He also provided a general
idea to define operations on soft sets. Maji et al. $^{\text{\cite{[2]}}}$
defined the operations of union and intersection on soft sets but Ali et al. 
$^{\text{\cite{[3]}}}$ pointed out some basic problems occurring with the
results using the definitions of \cite{[2]} and defined some new operations
on soft sets which were based upon the selection of parameters in
particular. Ali et al. $^{\text{\cite{[4]}}}$ studied the algebraic
structures of soft sets under these new operations. The hybrid structures
like fuzzy soft sets$^{\text{\cite{[5]}}}$, rough soft sets $^{\text{\cite%
{[6]}, \cite{[7]}}}$, intuitionistic fuzzy soft sets$^{\text{\cite{[8]}}}$,
vague soft sets $^{\text{\cite{[9]}}}$, and interval valued fuzzy soft sets $%
^{\text{\cite{[10]}}}$ were introduced and studied well through their
applications in various fields.

We have mentioned earlier that the motivation for the idea of soft sets is
the need of a parameterization tool and therefore, the role of parameters
becomes important and considerable. If $U$ is an initial universe and $E$ is
a set of parameter then a soft set defines a parameterized family of subsets
of $U$. Molodtsov presented an example in \cite{[1]} which discusses the
approximation of "attractiveness" of some houses through different
parameters like being expensive, beautiful, in the green surroundings etc.
Each approximation was taken as a set of houses possessing the property as
defined by the corresponding parameter. In that example the set of
parameters was $\{$expensive, beautiful, wooden, cheap, in the green
surroundings, modern, in good repair, in bad repair$\}$. We observe that the
parameters "expensive" and "cheap" and similarly "in good repair" and "in
bad repair" are clearly interrelated. The presence or absence of one affects
that of the other and so, a further filtration may be done to choose a
reduced set of parameters. An intelligent knowledge-based system must be
designed in a way that eliminates all the redundancies of input data and
hence it should provide a maximized and efficient level of performance. Each
parameter is a word or sentence and we need an efficient choice of
parameters as well. How to choose a precise, reduced and sufficient set of
parameters in order to optimize for a best possible decision? Chen et al. $^{%
\text{\cite{[11]}}}$ discussed a parameterization reduction of soft sets by
a step-wise reduction algorithm for parameters and compared it with the
concept of attribute reduction in rough sets. More recently M. I. Ali$^{%
\text{\cite{[12]}}}$ has given a new method for the reduction of parameters
in soft sets by keeping the equivalence of soft sets and information systems
in view. However, these discussions are mainly about the roles of decision
values i.e. the presence or absence of a particular parameter in a house and
did not consider the definition of parameter's set itself. Our question has
been somewhat answered by Dubois and Prade in \cite{[13]} while considering
the human judgments upon given pieces of information. They defined the role
of polarity by saying that choices of people are usually caused by checking
the good sides and the bad sides of alternatives separately. Then they
choose according to whether the good or the bad sides are stronger. In other
words the judgments also possess an intrinsic positive or negative flavor,
what we call a polarity. Bipolarity refers to an explicit handling of
positive and negative sides of information. Three types of bipolarity were
discussed in \cite{[13]} but we are using a rather generalized bipolarity
here, dealing with the positive and negative impacts only.

First of all we have to differentiate between a word or sentence represented
by a parameter and its negative (opposite) parameter e.g. \textquotedblleft
in good repair" and its opposite as \textquotedblleft in bad repair". In 
\cite{[2]} Maji et al. considered "not set" of parameters as $\{$not
expensive, not beautiful, not wooden, not cheap, not in the green
surroundings, not modern, not in good repair, not in bad repair$\}^{\text{%
\cite{[2]}}}$. It is worth noting that the parameters "not cheap", "not in a
good repair" and "expensive", "in a bad repair" are similar in meaning
respectively! This indicates the problem while defining the set of
parameters. In order to handle this problem we define a new concept of
bipolar soft set. A bipolar soft set is obtained by considering not only a
carefully chosen set of parameters but also an allied set of oppositely
meaning parameters named as "Not set of parameters". Structure of a bipolar
soft set is managed by two functions, say $F:A\rightarrow \mathcal{P}(U)$
and $G:\lnot A\rightarrow \mathcal{P}(U)$ where $\lnot A$ stands for the
"not set of $A$" and $G$ describes somewhat an opposite or negative
approximation for the attractiveness of a house relative to the
approximation computed by $F$. Maji et al.$^{\text{\cite{[2]}}}$ had used
the "not set" to define the compliment of a soft set. The \textit{complement
of a soft set }simply gives the complements of the approximations. The above
mentioned soft function $G$ is rather more generalized than soft complement
function and $(G$,$\lnot A)^{\text{\cite{[2]}}}$ can be any soft subset of $%
(F$,$A)^{c\text{ \cite{[2]}}}$. The difference is the gray area of choice,
that is, we may find some houses which do not satisfy any criteria properly
e.g. A house may not be highly expensive but it does not assure its
cheapness either. Thus, we must be careful while making our considerations
for the parameterization of data keeping in view that, during
approximations, there might be some indifferent elements in $U$.

We have given some preliminaries on soft sets and after that we define
bipolar soft sets. We have defined operations of union and intersection for
bipolar soft sets by taking restricted and extended sets of parameters. The
set of parameters has been taken to be restricted by taking intersection and
extended with the help of union and product operations. Examples are
presented which elaborate the concepts and working methods for computations.
The last section gives an application of bipolar soft sets in decision
making problems. We observe that the "Preference" is mostly viewed as an
important concept in decision making and the choice parameters are not all
equally attractive. Hence, it is crucial for a decision maker to examine
these parameters in terms of their desirability. Weights may be added to
parameters after comparing them on a common scale and thus we obtain a
weighted table for a bipolar soft set. A revised form of decision making
algorithm has also been presented and a comparison in results has been made
to clarify the use of weights.

\section{Soft Sets}

Let $U$ be an initial universe, and $E$ be a set of parameters. Let $%
\mathcal{P}(U)$ denotes the power set of $U$ and $A$, $B$ be non-empty
subsets of $E$.

\begin{definition}
\cite{[1]} A pair $(F,A)$\ is called a \textit{soft set} over $U$, where $F$
is a mapping given by $F:A\rightarrow \mathcal{P}(U)$.
\end{definition}

In other words, a soft set over $U$ is a parameterized family of subsets of
the universe $U$. For $e\in A$, $F(e)$ may be considered as the set of $e-$%
approximate elements of the soft set $(F$,$A)$. Clearly, a soft set is not a
set.

\begin{definition}
\cite{[3]} For two soft sets $(F,A)\ $and $(G,B)$ over a common universe $U$%
, we say that $(F,A)$ is a \textit{soft subset} of $(G,B)$ if

\begin{enumerate}
\item $A\subseteq B$ and

\item $F(e)\subseteq G(e)$ for all $e\in A$.
\end{enumerate}
\end{definition}

We write $(F,A)\widetilde{\subset }(G,B)$.

$(F,A)$ is said to be a \textit{soft super set} of $(G,B)$, if $(G,B)$ is a
soft subset of $(F,A)$. We denote it by $(F,A)\widetilde{\supset }(G,B)$.

\begin{definition}
\cite{[2]} Two soft sets $(F,A)$ and $(G,B)$ over a common universe $U$,\
are said to be \textit{soft equal} if $(F,A)$ is a soft subset of $(G,B)$\
and $(G,B)$ is a soft subset of $(F,A)$.
\end{definition}

\begin{definition}
$^{\text{\cite{[2]}}}$ Let\ $E=\{e_{1}$,$e_{2}$,...,$e_{n}\}$ be a set of
parameters. The $NOT$\textit{\ set of }$E$ denoted by$\ \lnot E$ is defined
by $\lnot E=\{\lnot e_{1}$,$\lnot e_{2}$,...,$\lnot e_{n}\}$ where, $\lnot
e_{i}=not$ $e_{i}$ for all $i$.
\end{definition}

\begin{proposition}
$^{\text{\cite{[2]}}}$\ For any subsets $A,B\subset E$,

\begin{enumerate}
\item $\lnot (\lnot A)=A$;

\item $\lnot (A\cup B)=\lnot A\cup \lnot B$;

\item $\lnot (A\cap B)=\lnot A\cap \lnot B$.
\end{enumerate}
\end{proposition}

\begin{definition}
$^{\text{\cite{[2]}}}$ The \textit{complement of a soft set }$(F$,$A)$ is
denoted by$\ (F$,$A)^{c}$\ and is defined by $(F$,$A)^{c}=(F^{c}$,$\lnot A)$
where, $F^{c}:\lnot A\rightarrow \mathcal{P}(U)$ is a mapping given by

$F^{c}(e)=U-F(\lnot e)$, for all $e\in \lnot A$.
\end{definition}

We call $F^{c}$ to be the soft complement function of $F$. Clearly $%
(F^{c})^{c}$ is the same as $F$ and $((F$,$A)^{c})^{c}=(F$,$A)$.

\begin{definition}
\cite{[3]} \textit{Union of two soft sets} $(F,A)$ and $(G,B)$ over the
common universe $U$\ is the soft set $(H,C)$, where $C=A\cup B$ and for all $%
e\in C$,%
\begin{equation*}
H(e)=\left\{ 
\begin{array}{cc}
F(e) & if\text{ }e\in A-B \\ 
G(e) & if\text{ }e\in B-A \\ 
F(e)\cup G(e) & if\text{ }e\in A\cap B%
\end{array}%
\right. 
\end{equation*}%
We write $(F,A)\cup _{\mathcal{E}}(G,B)=(H,C)$.
\end{definition}

\begin{definition}
\cite{[3]} \textit{The extended intersection} of two soft sets $(F,A)$ and $%
(G,B)$ over a common universe $U$, is the soft set $(H,C)$ where $C=\ A\cup B
$ and for all $e\in C$,%
\begin{equation*}
H\left( e\right) =\left\{ 
\begin{array}{cc}
F\left( e\right)  & if\ e\in A-B \\ 
G\left( e\right)  & if\ e\in B-A \\ 
F\left( e\right) \cap G\left( e\right)  & if\ e\in A\cap B%
\end{array}%
\right. .
\end{equation*}%
We write $(F,A)\cap _{\mathcal{E}}(G,B)=(H,C)$.
\end{definition}

\begin{definition}
\cite{[3]} Let $(F,A)$ and $(G,B)$ be two soft sets over the same universe $U
$\ such that $A\cap B\neq \emptyset $. The \textit{restricted intersection}
of $(F,A)$ and $(G,B)$ is denoted by $(F,A)\cap _{\mathcal{R}}(G,B)$ and is
defined as $(F,A)\cap _{\mathcal{R}}(G,B)=(H,A\cap B)$ where $H(e)=F(e)\cap
G(e)$ for all $e\in A\cap B$.
\end{definition}

\begin{definition}
\cite{[2]} Let $(F,A)$ and $(G,B)$ be two soft sets over the same universe $%
U $\ such that $A\cap B\neq \emptyset $. The \textit{restricted union} of $%
(F,A)$ and $(G,B)$ is denoted by $(F,A)\cup _{\mathcal{R}}(G,B)$ and is
defined as $(F,A)\cup _{\mathcal{R}}(G,B)=(H,C)$ where $C=A\cap B$ and for
all $e\in C$, $H(e)=F(e)\cup G(e)$.
\end{definition}

\begin{definition}
\cite{[3]} Let $U$ be an initial universe set, $E$ be the set of parameters,
and $A\subset E$.

\begin{enumerate}
\item $(F,A)$ is called a relative null soft set (with respect to the
parameter set $A$), denoted by $\Phi _{A}$, if $F(e)=\emptyset $ for all $%
e\in A$.

\item $(G,A)$ is called a relative whole soft set (with respect to the
parameter set $A$), denoted by $\mathfrak{U}_{A}$, if $F(e)=U$ for all $e\in
A$.
\end{enumerate}
\end{definition}

\section{Bipolar Soft Sets}

Let $U$ be an initial universe and $E$ be a set of parameters. Let $\mathcal{%
P}(U)$ denotes the power set of $U$ and $A$, $B$, $C$ are non-empty subsets
of $E$. For the convenience of reading, we shall refer the reader to \cite%
{[1]}, \cite{[2]} and \cite{[3]} for basic concepts of soft sets and their
properties. Now, we define

\begin{definition}
A triplet $(F$,$G$,$A)$ is called a \textit{bipolar} soft set over $U$,
where $F$ and $G$\ are mappings, given by $F:A\rightarrow \mathcal{P}(U)$
and\ $G:\lnot A\rightarrow \mathcal{P}(U)$ such that $F(e)\cap G(\lnot
e)=\emptyset $ (Empty Set) for all $e\in A$.
\end{definition}

In other words, a \textit{bipolar} soft set over $U$\ gives two parametrized
families of subsets of the universe $U$ and the condition $F(e)\cap G(\lnot
e)=\emptyset $ for all $e\in A$, is imposed as a consistency constraint. For
each $e\in A$, $F(e)$ and $G(\lnot e)$ are regarded as the set of $e$%
-approximate elements of the \textit{bipolar} soft set $(F$,$G$,$A)$. It is
also observed that the relationship between the complement function and the
defining function of a soft set is a particular case for the defining
functions of a bipolar soft set, i.e. $(F$,$F^{c}$,$A)$ is a bipolar soft
set over $U$. The difference occurs due to the presence of uncertainty or
hesitation or the lack of knowledge in defining the membership function. We
name this uncertainty or gray area as the approximation for the degree of
hesitation. Thus the union of three approximations i.e. $e$-approximation, $%
\lnot e$-approximation, and approximation of hesitation is $U$. We note that 
$\emptyset \subseteq U-\{F(e)\cup G(\lnot e)\}\subseteq U$, for each $e\in A$%
. So, we may approximate the degree of hesitation in $(F$,$G$,$A)$ by an
allied soft set $(H$,$A)$ defined over $U$, where $H(e)=U-\{F(e)\cup G(\lnot
e)\}$ for all $e\in A$.

\begin{definition}
For two \textit{bipolar} soft sets $(F$,$G$,$A)$ and $(F_{1}$,$G_{1}$,$B)$
over a universe $U$, we say that $(F$,$G$,$A)$ is a \textit{bipolar} soft
subset of $(F_{1}$,$G_{1}$,$B)$, if,

\begin{enumerate}
\item $A\subseteq B$ and

\item $F(e)\subseteq F_{1}(e)$ and $G_{1}(\lnot e)\subseteq G(\lnot e)$ for
all $e\in A$.
\end{enumerate}
\end{definition}

This relationship is denoted by $(F$,$G$,$A)\tilde{\subseteq}(F_{1}$,$G_{1}$,%
$B)$. Similarly $(F$,$G$,$A)$ is said to be a \textit{bipolar} soft superset
of $(F_{1}$,$G_{1}$,$B)$, if $(F_{1,}G_{1}$,$B)$ is a \textit{bipolar soft
subset} of $(F$,$G$,$A)$. We denote it by $(F$,$G$,$A)\tilde{\supseteq}%
(F_{1,}G_{1}$,$B)$.

\begin{definition}
Two \textit{bipolar} soft sets $(F$,$G$,$A)$ and $(F_{1}$,$G_{1}$,$B)$ over
a universe $U$ are said to be \textit{equal} if $(F$,$G$,$A)$ is a \textit{%
bipolar} soft subset of $(F_{1}$,$G_{1}$,$B)$ and $(F_{1}$,$G_{1}$,$B)$ is a 
\textit{bipolar} soft subset of $(F$,$G$,$A)$.
\end{definition}

\begin{definition}
The \textit{complement of a bipolar soft set} $(F$,$G$,$A)$ is denoted by $%
(F $,$G$,$A)^{c}$ and is defined by $(F$,$G$,$A)^{c}=(F^{c}$,$G^{c}$,$A)$
where $F^{c}$ and $G^{c}$\ are mappings given by $F^{c}(e)=G(\lnot e)$ and $%
G^{c}(\lnot e)=F(e)$ for all $e\in A$.
\end{definition}

\begin{definition}
A \textit{bipolar} soft set over $U$ is said to be a \textit{relative null
bipolar soft set,} denoted by $(\Phi $,$\mathfrak{U}$,$A)$ if for all $e\in
A $, $\Phi (e)=\emptyset $ and $\mathfrak{U}(\lnot e)=U$, for all $e\in A$.
\end{definition}

\begin{definition}
A \textit{bipolar} soft set over $U$ is said to be a \textit{relative
absolute bipolar soft set,} denoted by $(\mathfrak{U}$,$\Phi $,$A)$, if for
all $e\in A$, $\mathfrak{U}(e)=U$ and $\Phi (\lnot e)=\emptyset $, for all $%
e\in A$.
\end{definition}

\begin{definition}
If $(F$,$G$,$A)$ and $(F_{1}$,$G_{1}$,$B)$ are two \textit{bipolar} soft
sets over $U$ then " $(F$,$G$,$A)$ and $(F_{1}$,$G_{1}$,$B)$" denoted by $(F$%
,$G$,$A)\wedge (F_{1}$,$G_{1}$,$B)$ is defined by%
\begin{equation*}
(F,G,A)\wedge (F_{1},G_{1},B)=(H,I,A\times B)
\end{equation*}%
where $H(a$,$b)=F(a)\cap F_{1}(b)$ and $I(\lnot a$,$\lnot b)=G(\lnot a)\cup
G_{1}(\lnot b)$,\ for all $(a$,$b)\in A\times B$.
\end{definition}

\begin{definition}
If $(F$,$G$,$A)$ and $(F_{1}$,$G_{1}$,$B)$ are two \textit{bipolar} soft
sets over $U$ then " $(F$,$G$,$A)$ or $(F_{1}$,$G_{1}$,$B)$" denoted by $(F$,%
$G$,$A)\vee (F_{1}$,$G_{1}$,$B)$ is defined by%
\begin{equation*}
(F,G,A)\vee (F_{1},G_{1},B)=(H,I,A\times B)
\end{equation*}%
where $H(a$,$b)=F(a)\cup F_{1}(b)$ and $I(\lnot a$,$\lnot b)=G(\lnot a)\cap
G_{1}(\lnot b)$,\ for all $(a$,$b)\in A\times B$.
\end{definition}

\begin{proposition}
If $(F$,$G$,$A)$ and $(F_{1}$,$G_{1}$,$B)$\ are two \textit{bipolar} soft
sets over $U$ then

\begin{enumerate}
\item $((F$,$G$,$A)\vee (F_{1}$,$G_{1}$,$B))^{c}=(F$,$G$,$A)^{c}\wedge
(F_{1} $,$G_{1}$,$B)^{c}$

\item $((F$,$G$,$A)\wedge (F_{1}$,$G_{1}$,$B))^{c}=(F$,$G$,$A)^{c}\vee
(F_{1} $,$G_{1}$,$B)^{c}$.
\end{enumerate}
\end{proposition}

\begin{proof}
\ Straightforward.\ 
\end{proof}

\begin{definition}
\textit{Extended Union of two bipolar soft sets} $(F$,$G$,$A)$ and $(F_{1}$,$%
G_{1}$,$B)$ over the common universe $U$\ is the \textit{bipolar} soft set $%
(H$,$I$,$C)$ over $U$, where $C=A\cup B$ and for all $e\in C$,%
\begin{eqnarray*}
H(e) &=&\left\{ 
\begin{array}{cc}
F(e) & if\text{ }e\in A-B \\ 
F_{1}(e) & if\text{ }e\in B-A \\ 
F(e)\cup F_{1}(e) & if\text{ }e\in A\cap B%
\end{array}%
\right. \\
I(\lnot e) &=&\left\{ 
\begin{array}{cc}
G(\lnot e) & if\text{ }e\in (\lnot A)-(\lnot B) \\ 
G_{1}(\lnot e) & if\text{ }e\in (\lnot B)-(\lnot A) \\ 
G(\lnot e)\cap G_{1}(\lnot e) & if\text{ }e\in (\lnot A)\cap (\lnot B)%
\end{array}%
\right.
\end{eqnarray*}%
We denote it by $(F$,$G$,$A)\tilde{\cup}(F_{1}$,$G_{1}$,$B)=(H$,$I$,$C)$.
\end{definition}

\begin{definition}
\textit{Extended Intersection of two bipolar soft sets} $(F$,$G$,$A)$ and $%
(F_{1}$,$G_{1}$,$B)$ over the common universe $U$\ is the \textit{bipolar}
soft set $(H$,$I$,$C)$ over $U$, where $C=A\cup B$ and for all $e\in C$,%
\begin{eqnarray*}
H(e) &=&\left\{ 
\begin{array}{cc}
F(e) & if\text{ }e\in A-B \\ 
F_{1}(e) & if\text{ }e\in B-A \\ 
F(e)\cap F_{1}(e) & if\text{ }e\in A\cap B%
\end{array}%
\right. \\
I(\lnot e) &=&\left\{ 
\begin{array}{cc}
G(\lnot e) & if\text{ }e\in (\lnot A)-(\lnot B) \\ 
G_{1}(\lnot e) & if\text{ }e\in (\lnot B)-(\lnot A) \\ 
G(\lnot e)\cup G_{1}(\lnot e) & if\text{ }e\in (\lnot A)\cap (\lnot B)%
\end{array}%
\right.
\end{eqnarray*}%
We denote it by $(F$,$G$,$A)\tilde{\cap}(F_{1}$,$G_{1}$,$B)=(H$,$I$,$C)$.
\end{definition}

\begin{definition}
\textit{Restricted Union of two bipolar soft sets} $(F$,$G$,$A)$ and $(F_{1}$%
,$G_{1}$,$B)$ over the common universe $U$\ is the \textit{bipolar} soft set 
$(H$,$I$,$C)$, where $C=A\cap B$ is non-empty and for all $e\in C$%
\begin{equation*}
H(e)=F(e)\cup G(e)\text{ and }I(\lnot e)=F_{1}(\lnot e)\cap G_{1}(\lnot e)%
\text{.}
\end{equation*}%
We denote it by $(F$,$G$,$A)\cup _{\mathcal{R}}(F_{1}$,$G_{1}$,$B)=(H$,$I$,$%
C)$.
\end{definition}

\begin{definition}
\textit{Restricted Intersection of two bipolar soft sets} $(F$,$G$,$A)$ and $%
(F_{1}$,$G_{1}$,$B)$ over the common universe $U$\ is the \textit{bipolar}
soft set $(H$,$I$,$C)$, where $C=A\cap B$ is non-empty and for all $e\in C$:%
\begin{equation*}
H(e)=F(e)\cap G(e)\text{ and }I(\lnot e)=F_{1}(\lnot e)\cup G_{1}(\lnot e)%
\text{.}
\end{equation*}%
We denote it by $(F$,$G$,$A)\cap _{\mathcal{R}}(F_{1}$,$G_{1}$,$B)=(H$,$I$,$%
C)$.
\end{definition}

\begin{proposition}
Let $(F$, $G$, $A)$ and $(F_{1}$,$G_{1}$,$A)$\ be two \textit{bipolar} soft
sets over a common universe $U$. Then the following are true

\begin{enumerate}
\item $((F$,$G$,$A)\tilde{\cup}(F_{1}$,$G_{1}$,$B))^{c}=(F$,$G$,$A)^{c}%
\tilde{\cap}(F_{1}$,$G_{1}$,$B)^{c}$,

\item $((F$,$G$,$A)\tilde{\cap}(F_{1}$,$G_{1}$,$B))^{c}=(F$,$G$,$A)^{c}%
\tilde{\cup}(F_{1}$,$G_{1}$,$B)^{c}$,

\item $((F$,$G$,$A)\cup _{\mathcal{R}}(F_{1}$,$G_{1}$,$B))^{c}=(F$,$G$,$%
A)^{c}\cap _{\mathcal{R}}(F_{1}$,$G_{1}$,$B)^{c}$,

\item $((F$,$G$,$A)\cap _{\mathcal{R}}(F_{1}$,$G_{1}$,$B))^{c}=(F$,$G$,$%
A)^{c}\cup _{\mathcal{R}}(F_{1}$,$G_{1}$,$B)^{c}$.
\end{enumerate}

\begin{proof}
1) Let $e\in A\cup B$. There are three cases:

\begin{itemize}
\item[(i)] If $e\in A-B$, then 
\begin{eqnarray*}
(F\tilde{\cup}F_{1})^{c}(e) &=&(F(e))^{c}=(F^{c}\tilde{\cap}F_{1}{}^{c})(e)
\\
(G\tilde{\cup}G_{1})^{c}(e) &=&(G(e))^{c}=(G^{c}\tilde{\cap}G_{1}{}^{c})(e)%
\text{,}
\end{eqnarray*}

\item[(ii)] If $e\in B-A$, then 
\begin{eqnarray*}
(F\tilde{\cup}F_{1})^{c}(e) &=&(F_{1}(e))^{c}=(F^{c}\tilde{\cap}%
F_{1}{}^{c})(e) \\
(G\tilde{\cup}G_{1})^{c}(e) &=&(G_{1}(e))^{c}=(G^{c}\tilde{\cap}%
G_{1}{}^{c})(e)\text{,}
\end{eqnarray*}

\item[(iii)] If $e\in A\cap B$, then 
\begin{eqnarray*}
(F\tilde{\cup}F_{1})^{c}(e) &=&(F(e)\cup F_{1}(e))^{c}=(F(e))^{c}\cap
(F_{1}(e))^{c} \\
(G\tilde{\cup}G_{1})^{c}(e) &=&(G(e)\cap G_{1}(e))^{c}=(G(e))^{c}\cup
(G_{1}(e))^{c}\text{,}
\end{eqnarray*}%
and,%
\begin{eqnarray*}
(F^{c}\tilde{\cap}F_{1}{}^{c})(e) &=&(F(e))^{c}\cap (F_{1}(e))^{c} \\
(G^{c}\tilde{\cap}G_{1}{}^{c})(e) &=&(G(e))^{c}\cup (G_{1}(e))^{c}\text{,}
\end{eqnarray*}
\end{itemize}

Therefore, in all three cases we obtain equality and thus 
\begin{equation*}
((F,G,A)\tilde{\cup}(F_{1},G_{1},B))^{c}=(F,G,A)^{c}\tilde{\cap}%
(F_{1},G_{1},B)^{c}\text{.}
\end{equation*}%
The remaining parts can also be proved in a similar way.
\end{proof}
\end{proposition}

\begin{proposition}
If $(\Phi $,$\mathfrak{U}$,$A)$ is a null \textit{bipolar} soft set, $(%
\mathfrak{U}$,$\Phi $,$A)$ an absolute \textit{bipolar} soft set, and $(F$,$%
G $,$A)$, $(F_{1}$,$G_{1}$,$A)$ are \textit{bipolar} soft sets over $U$, then

\begin{enumerate}
\item $(F$,$G$,$A)\tilde{\cup}(F_{1}$,$G_{1}$,$A)=(F$,$G$,$A)\cup _{\mathcal{%
R}}(F_{1}$,$G_{1}$,$A)$,

\item $(F$,$G$,$A)\tilde{\cap}(F_{1}$,$G_{1}$,$A)=(F$,$G$,$A)\cap _{\mathcal{%
R}}(F_{1}$,$G_{1}$,$A)$,

\item $(F$,$G$,$A)\tilde{\cup}(F$,$G$,$A)=(F$,$G$,$A)$; $(F$,$G$,$A)\tilde{%
\cap}(F$,$G$,$A)=(F$,$G$,$A)$,

\item $(F$,$G$,$A)\tilde{\cup}(\Phi $,$\mathfrak{U}$,$A)=(F$,$G$,$A)$; $(F$,$%
G$,$A)\tilde{\cap}(\Phi $,$\mathfrak{U}$,$A)=(\Phi $,$\mathfrak{U}$,$A)$,

\item $(F$,$G$,$A)\tilde{\cup}(\mathfrak{U}$,$\Phi $,$A)=(\mathfrak{U}$,$%
\Phi $,$A)$; $(F$,$G$,$A)\tilde{\cap}(\mathfrak{U}$,$\Phi $,$A)=(F$,$G$,$A)$.
\end{enumerate}

\begin{proof}
\ Straightforward.\ 
\end{proof}
\end{proposition}

\begin{example}
Let $U$ be the set of houses under consideration, and $E$ be the set of
parameters, $U=\{h_{1}$,$h_{2}$,$h_{3}$,$h_{4}$,$h_{5}\}$,

$E=\{e_{1}$,$e_{2}$,$e_{3}$,$e_{4}$,$e_{5}$,$e_{6}\}=\{$ in the green
surroundings, wooden, cheap, in good repair, furnished, traditional $\}$.
Then $\lnot E=\{\lnot e_{1}$,$\lnot e_{2}$,$\lnot e_{3}$,$\lnot e_{4}$,$%
\lnot e_{5}$,$\lnot e_{6}\}=\{$ in the commercial area, marbled, expensive,
in bad repair, non-furnished, modern $\}$. Suppose that $A=\{e_{1}$,$e_{2}$,$%
e_{3}$,$e_{6}\}$, and $B=\{e_{2}$,$e_{3}$,$e_{4}$,$e_{5}\}$. The bipolar
soft sets $(F$,$G$,$A)$ and $(F_{1}$,$G_{1}$,$B)$ describe the
\textquotedblleft requirements of the houses\textquotedblright\ which Mr. X
and Mr. Y are going to buy respectively. Suppose that%
\begin{eqnarray*}
F(e_{1})
&=&\{h_{2},h_{3}\},F(e_{2})=\{h_{1},h_{2},h_{5}\},F(e_{3})=\{h_{1},h_{3}%
\},F(e_{6})=\{h_{2},h_{3},h_{5}\} \\
G(\lnot e_{1}) &=&\{h_{4},h_{5}\},G(\lnot e_{2})=\{h_{3},h_{4}\},G(\lnot
e_{3})=\{h_{2},h_{4}\},G(\lnot e_{6})=\{h_{4}\}
\end{eqnarray*}%
and%
\begin{align*}
F_{1}(e_{2})&
=\{h_{2},h_{5}\},F_{1}(e_{3})=\{h_{1},h_{3},h_{5}\},F_{1}(e_{4})=%
\{h_{1},h_{3},h_{4}\},F_{1}(e_{5})=\{h_{2},h_{3}\}, \\
G_{1}(\lnot e_{2})& =\{h_{4}\},G_{1}(\lnot
e_{3})=\{h_{2},h_{4}\},G_{1}(\lnot e_{4})=\{h_{2}\},G_{1}(\lnot
e_{5})=\{h_{1},h_{4}\}.
\end{align*}%
One may read this information as: Mr. X demands a house according to the
parameters, " in the green surroundings, wooden, cheap and traditional". The
houses in the set $U$ are parametrized accordingly, either having the said
attribute or its opposite. We observe that Mr. X thinks that the houses $%
h_{2}$ and $h_{3}$ are situated in the green surroundings while $h_{4}$ and $%
h_{5}$ are situated in an entirely commercial area and we are not
parameterizing $h_{1}$\ by either of $A$ or $\lnot A$. It means that $h_{1}$
has some green surrounding area and some area which is either commercial or
not green as such. This also shows that the approximations for $(G$,$\lnot A)
$ are not same as that of $(F$,$A)^{c}$ as defined in \cite{[2]}. We note
that, for Mr. Y, the houses $h_{2}$ and $h_{5}$ are wooden while $h_{4}$ is
marbled whereas $h_{1}$ and $h_{3}$ are not sufficiently structured to be
parametrized as either of these. We also see that Mr. X thinks that $h_{1}$
and $h_{3}$ are cheap while $h_{2}$ and $h_{4}$ are expensive but Mr. Y
considers that $h_{1}$, $h_{3}$ and $h_{4}$ are cheap while $h_{2}$ and $%
h_{4}$ are expensive ones. So parameterization differs according to any
given situation and circumstances and further approximations about the given
data may be obtained by applying the operations of union and intersection
defined for bipolar soft sets.

Now, we approximate the resulting bipolar soft sets obtained by applying the
above mentioned operations on $(F$,$G$,$A)$ and $(F_{1}$,$G_{1}$,$B)$.

Let $(F,G,A)\tilde{\cup}(F_{1},G_{1},B)=(H_{1},I_{1},A\cup B)$. Then%
\begin{eqnarray*}
H_{1}(e_{1})
&=&\{h_{2},h_{3}\},H_{1}(e_{2})=\{h_{1},h_{2},h_{5}\},H_{1}(e_{3})=%
\{h_{1},h_{3},h_{5}\}, \\
H_{1}(e_{4})
&=&\{h_{1},h_{3},h_{4}\},H_{1}(e_{5})=\{h_{2},h_{3}\},H_{1}(e_{6})=%
\{h_{2},h_{3},h_{5}\},
\end{eqnarray*}%
and%
\begin{eqnarray*}
I_{1}(\lnot e_{1}) &=&\{h_{4},h_{5}\},I_{1}(\lnot
e_{2})=\{h_{4}\},I_{1}(\lnot e_{3})=\{h_{2},h_{4}\},I_{1}(\lnot
e_{4})=\{h_{4}\}, \\
I_{1}(\lnot e_{5}) &=&\{h_{1},h_{4}\},I_{1}(\lnot e_{6})=\{h_{4}\}.
\end{eqnarray*}%
Let $(F,G,A)\tilde{\cap}(F_{1},G_{1},B)=(H_{2},I_{2},A\cup B)$. Then%
\begin{eqnarray*}
H_{2}(e_{1})
&=&\{h_{2},h_{3}\},H_{2}(e_{2})=\{h_{2},h_{5}\},H_{2}(e_{3})=\{h_{1},h_{3}%
\},H_{2}(e_{4})=\{h_{1},h_{3},h_{4}\}, \\
H_{2}(e_{5}) &=&\{h_{2},h_{3}\},H_{2}(e_{6})=\{h_{2},h_{3},h_{5}\},
\end{eqnarray*}%
and%
\begin{eqnarray*}
I_{2}(\lnot e_{1}) &=&\{h_{4},h_{5}\},I_{2}(\lnot
e_{2})=\{h_{3},h_{4}\},I_{2}(\lnot e_{3})=\{h_{2},h_{4}\},I_{2}(\lnot
e_{4})=\{h_{2}\}, \\
I_{2}(\lnot e_{5}) &=&\{h_{1},h_{4}\},I_{2}(\lnot e_{6})=\{h_{4}\}.
\end{eqnarray*}%
Let $(F,G,A)\cup _{\mathcal{R}}(F_{1},G_{1},B)=(H_{3},I_{3},A\cap B)$. Then%
\begin{eqnarray*}
H_{3}(e_{2}) &=&\{h_{1},h_{2},h_{5}\},H_{3}(e_{3})=\{h_{1},h_{3},h_{5}\}%
\text{ \ \ \ and} \\
I_{3}(\lnot e_{2}) &=&\{h_{4}\},I_{3}(\lnot e_{3})=\{h_{2},h_{4}\}.
\end{eqnarray*}%
Let $(F,G,A)\cap _{\mathcal{R}}(F_{1},G_{1},B)=(H_{4},I_{4},A\cap B)$. Then%
\begin{eqnarray*}
H_{4}(e_{2}) &=&\{h_{2},h_{5}\},H_{4}(e_{3})=\{h_{1},h_{3}\}\text{ \ and} \\
I_{4}(\lnot e_{2}) &=&\{h_{3},h_{4}\},I_{4}(\lnot e_{3})=\{h_{2},h_{4}\}.
\end{eqnarray*}%
Let $(F,G,A)\vee (F_{1},G_{1},B)=(H_{5},I_{5},A\times B)$. Then%
\begin{eqnarray*}
H_{5}(e_{1},e_{2})
&=&\{h_{2},h_{3},h_{5}\},H_{5}(e_{1},e_{3})=\{h_{1},h_{2},h_{3},h_{5}%
\},H_{5}(e_{1},e_{4})=\{h_{1},h_{2},h_{3},h_{4}\}, \\
H_{5}(e_{1},e_{5})
&=&\{h_{2},h_{3}\},H_{5}(e_{2},e_{2})=\{h_{1},h_{2},h_{5}%
\},H_{5}(e_{2},e_{3})=\{h_{1},h_{2},h_{3},h_{5}\}
\end{eqnarray*}
and%
\begin{eqnarray*}
I_{5}(\lnot e_{1},\lnot e_{2}) &=&\{h_{4}\},I_{5}(\lnot e_{1},\lnot
e_{3})=\{h_{4}\},I_{5}(\lnot e_{1},\lnot e_{4})=\emptyset \text{,}%
I_{5}(\lnot e_{1},\lnot e_{5})=\{h_{4}\}, \\
I_{5}(\lnot e_{2},\lnot e_{2}) &=&\{h_{4}\},I_{5}(\lnot e_{2},\lnot
e_{3})=\{h_{4}\}\text{ \ and so on.}
\end{eqnarray*}%
Let $(F,G,A)\wedge (F_{1},G_{1},B)=(H_{6},I_{6},A\times B)$. Then%
\begin{eqnarray*}
H_{6}(e_{1},e_{2})
&=&\{h_{2}\},H_{6}(e_{1},e_{3})=\{h_{3}\},H_{6}(e_{1},e_{4})=\{h_{3}%
\},H_{6}(e_{1},e_{5})=\{h_{2},h_{3}\}, \\
H_{6}(e_{2},e_{2}) &=&\{h_{2},h_{5}\},H_{6}(e_{2},e_{3})=\{h_{1},h_{5}\}
\end{eqnarray*}
and%
\begin{eqnarray*}
I_{6}(\lnot e_{1},\lnot e_{2}) &=&\{h_{4},h_{5}\},I_{6}(\lnot e_{1},\lnot
e_{3})=\{h_{2},h_{4},h_{5}\},I_{6}(\lnot e_{1},\lnot
e_{4})=\{h_{2},h_{4},h_{5}\}, \\
I_{6}(\lnot e_{1},\lnot e_{5}) &=&\{h_{1},h_{4},h_{5}\},I_{6}(\lnot
e_{2},\lnot e_{2})=\{h_{3},h_{4}\},I_{6}(\lnot e_{2},\lnot
e_{3})=\{h_{2},h_{3},h_{4}\}
\end{eqnarray*}%
and so on.
\end{example}

\label{tbldesc1}A bipolar soft set may be represented by a pair of tables
for each of the functions $F$ and $G$ respectively in a similar way as the
tabular representation of soft sets is used by Maji et al. in \cite{[2]}. We
can also represent a bipolar soft set with the help of a single table by
putting%
\begin{equation*}
a_{ij}=\left\{ 
\begin{array}{rll}
1 &  & \text{if }h_{i}\in F(e_{j}) \\ 
0 &  & \text{if }h_{i}\in U-\{F(e_{j})\cup G(\lnot e_{j})\} \\ 
-1 &  & \text{if }h_{i}\in G(\lnot e_{j})%
\end{array}%
\right. 
\end{equation*}%
where $a_{ij}$ is the $ith$ entry of $jth$ column of the table whose rows
and columns are labeled by houses and parameters respectively. The tabular
representations of the bipolar soft set $(F,G,A)$ are given by Table \ref%
{TableKey copy(1)} and Table \ref{TableKey copy(2)}.

\begin{table}[h] \centering%
$%
\begin{tabular}{|l|l|l|l|l|}
\hline
\textbf{F} & $\mathbf{e}_{1}$ & $\mathbf{e}_{2}$ & $\mathbf{e}_{3}$ & $%
\mathbf{e}_{6}$ \\ \hline
$h_{1}$ & $0$ & $1$ & $1$ & $0$ \\ \hline
$h_{2}$ & $1$ & $1$ & $0$ & $1$ \\ \hline
$h_{3}$ & $1$ & $0$ & $1$ & $1$ \\ \hline
$h_{4}$ & $0$ & $0$ & $0$ & $0$ \\ \hline
$h_{5}$ & $0$ & $1$ & $0$ & $1$ \\ \hline
\end{tabular}%
$ $\ \ \ \ \ \ \ \ 
\begin{tabular}{|l|l|l|l|l|}
\hline
\textbf{G} & $\lnot \mathbf{e}_{1}$ & $\lnot \mathbf{e}_{2}$ & $\lnot 
\mathbf{e}_{3}$ & $\lnot \mathbf{e}_{6}$ \\ \hline
$h_{1}$ & $0$ & $0$ & $0$ & $0$ \\ \hline
$h_{2}$ & $0$ & $0$ & $1$ & $0$ \\ \hline
$h_{3}$ & $0$ & $1$ & $0$ & $0$ \\ \hline
$h_{4}$ & $1$ & $1$ & $1$ & $1$ \\ \hline
$h_{5}$ & $1$ & $0$ & $0$ & $0$ \\ \hline
\end{tabular}%
$\caption{Tabular Representaion of (F,G,A) using a Pair of Tables}\label%
{TableKey copy(1)}%
\end{table}
\begin{table}[h] \centering%
$%
\begin{tabular}{|l|l|l|l|l|}
\hline
(\textbf{F,G,A)} & $\mathbf{e}_{1}$ & $\mathbf{e}_{2}$ & $\mathbf{e}_{3}$ & $%
\mathbf{e}_{6}$ \\ \hline
$h_{1}$ & \multicolumn{1}{|r|}{$0$} & \multicolumn{1}{|r|}{$1$} & 
\multicolumn{1}{|r|}{$1$} & \multicolumn{1}{|r|}{$0$} \\ \hline
$h_{2}$ & \multicolumn{1}{|r|}{$1$} & \multicolumn{1}{|r|}{$1$} & 
\multicolumn{1}{|r|}{$-1$} & \multicolumn{1}{|r|}{$1$} \\ \hline
$h_{3}$ & \multicolumn{1}{|r|}{$1$} & \multicolumn{1}{|r|}{$-1$} & 
\multicolumn{1}{|r|}{$1$} & \multicolumn{1}{|r|}{$1$} \\ \hline
$h_{4}$ & \multicolumn{1}{|r|}{$-1$} & \multicolumn{1}{|r|}{$-1$} & 
\multicolumn{1}{|r|}{$-1$} & \multicolumn{1}{|r|}{$-1$} \\ \hline
$h_{5}$ & \multicolumn{1}{|r|}{$-1$} & \multicolumn{1}{|r|}{$1$} & 
\multicolumn{1}{|r|}{$0$} & \multicolumn{1}{|r|}{$1$} \\ \hline
\end{tabular}%
$\caption{Tabular Representaion of (F,G,A) using only one Table}\label%
{TableKey copy(2)}%
\end{table}%

\begin{example}
Bipolar disorder is a serious psychological illness that can lead to
dangerous behavior, problematic careers and relationships, and suicidal
tendencies, especially if not treated early. A bipolar mood chart is a
simple and yet effective means of tracking and representing patient's
condition every month. Bipolar mood charts help patients, their families and
their doctors to see probable patterns that might have been very difficult
to determine. Bipolar children and their families will greatly benefit from
mood charting and can expect early detection of symptoms and determination
of proper treatments by their doctors. We construct a mood chart based upon
a bipolar soft set as follows:

Let $U=\{1,2,3,4,5,6,7\}$ be the set of days in which the record has been
maintained i.e. $i=ith$ day under observation, for $1\leq i\leq 7$. Let

$E=\{e_{1},e_{2},e_{3},e_{4},e_{5}\}=\{$Severe Mania, Severe Depression,
Anxiety, Medication, Side effects$\}$ and

$\lnot E=\{\lnot e_{1},\lnot e_{2},\lnot e_{3},\lnot e_{4},\lnot e_{5}\}=\{$%
Mild Mania, Mild Depression, No Anxiety, No Medication, No Side effects$\}$.
Here the gray area is obviously the moderate form of parameters. Let the
bipolar soft sets $(F$,$G$,$E)$ describes the \textquotedblleft daily record
of the behavior\textquotedblright\ of Mr. X. Suppose that%
\begin{eqnarray*}
F(e_{1}) &=&\{1,5\}\text{, }F(e_{2})=\{1,2,3,4,7\}\text{, }%
F(e_{3})=\{2,4,5,6\}\text{,} \\
F(e_{4}) &=&\{1,2,4,5,6,7\}\text{, }F(e_{5})=\{2,3,5,7\}\text{.}
\end{eqnarray*}%
and%
\begin{eqnarray*}
G(\lnot e_{1}) &=&\{2,6,7\}\text{, }G(\lnot e_{2})=\{6\}\text{, }G(\lnot
e_{3})=\{1,7\}\text{, }G(\lnot e_{4})=\{3\}\text{, } \\
G(\lnot e_{5}) &=&\{1,4,6\}\text{.}
\end{eqnarray*}%
The corresponding mood chart is given by the Table \ref{TableKey}

\begin{table}[h] \centering%
\begin{tabular}[t]{|r|r|r|r|r|r|}
\hline
$\mathbf{(F,G,E)}$ & $e_{1}$ & $e_{2}$ & $e_{3}$ & $e_{4}$ & $e_{5}$ \\ 
\hline
$1$ & $1$ & $1$ & $-1$ & $1$ & $-1$ \\ \hline
$2$ & $-1$ & $1$ & $1$ & $1$ & $1$ \\ \hline
$3$ & $0$ & $1$ & $0$ & $-1$ & $1$ \\ \hline
$4$ & $0$ & $1$ & $1$ & $1$ & $-1$ \\ \hline
$5$ & $1$ & $0$ & $1$ & $1$ & $1$ \\ \hline
$6$ & $-1$ & $-1$ & $1$ & $1$ & $-1$ \\ \hline
$7$ & $-1$ & $1$ & $-1$ & $1$ & $1$ \\ \hline
\end{tabular}%
\caption{Tabular Representation of (F,G,E)}\label{TableKey}%
\end{table}%
\end{example}

\begin{lemma}
Let $(F,$ $G,$ $A)$, $(F_{1},G_{1},B)$ and $(F_{2},G_{2},C)$\ be any \textit{%
bipolar soft set}s over a common universe $U$. Then the following are true

\begin{enumerate}
\item $(F,G,A)\alpha ((F_{1},G_{1},B)\alpha (F_{2},G_{2},C))=((F,G,A)\alpha
(F_{1},G_{1},B))\alpha (F_{2},G_{2},C)$,

\item $(F,G,A)\alpha (F_{1},G_{1},B)=(F,G,A)\alpha (F_{1},G_{1},B)$,
\end{enumerate}

for all $\alpha \in \{\tilde{\cap},\cap _{\mathcal{R}},\tilde{\cup},\cup _{%
\mathcal{R}}\}$.

\begin{proof}
Straightforward.\ 
\end{proof}
\end{lemma}

\begin{lemma}
Let $(F,$ $G,$ $A)$ and $(F_{1},G_{1},B)$\ be two \textit{bipolar soft set}s
over a common universe $U$. Then the following are true

\begin{enumerate}
\item $(F,G,A)\tilde{\cup}(F_{1},G_{1},B)$ is the smallest \textit{bipolar}
soft set over $U$ which contains both $(F,G,A)$ and $(F_{1},G_{1},B)$.

\item $(F,G,A)\cap _{\mathcal{R}}(F_{1},G_{1},B)$ is the largest \textit{%
bipolar} soft set over $U$ which is contained in both $(F,G,A)$ and $%
(F_{1},G_{1},B)$.
\end{enumerate}

\begin{proof}
Straightforward.
\end{proof}
\end{lemma}

\begin{proposition}
Let $(F,$ $G,$ $A)$, $(F_{1},G_{1},B)$ and $(F_{2},G_{2},C)$\ be any \textit{%
bipolar soft set}s over a common universe $U$. Then

\begin{enumerate}
\item $(F,G,A)\cap _{\mathcal{R}}((F_{1},G_{1},B)\tilde{\cup}%
(F_{2},G_{2},C))=((F,G,A)\cap _{\mathcal{R}}(F_{1},G_{1},B))\tilde{\cup}%
((F,G,A)\cap _{\mathcal{R}}(F_{2},G_{2},C))$,

\item $(F,G,A)\cap _{\mathcal{R}}((F_{1},G_{1},B)\cap _{\mathcal{R}%
}(F_{2},G_{2},C))=((F,G,A)\cap _{\mathcal{R}}(F_{1},G_{1},B))\cap _{\mathcal{%
R}}((F,G,A)\cap _{\mathcal{R}}(F_{2},G_{2},C))$,

\item $(F,G,A)\cap _{\mathcal{R}}((F_{1},G_{1},B)\cup _{\mathcal{R}%
}(F_{2},G_{2},C))=((F,G,A)\cap _{\mathcal{R}}(F_{1},G_{1},B))\cup _{\mathcal{%
R}}((F,G,A)\cap _{\mathcal{R}}(F_{2},G_{2},C))$,

\item $(F,G,A)\cup _{\mathcal{R}}((F_{1},G_{1},B)\tilde{\cup}%
(F_{2},G_{2},C))=((F,G,A)\cup _{\mathcal{R}}(F_{1},G_{1},B))\tilde{\cup}%
((F,G,A)\cup _{\mathcal{R}}(F_{2},G_{2},C))$,

\item $(F,G,A)\cup _{\mathcal{R}}((F_{1},G_{1},B)\cap _{\mathcal{R}%
}(F_{2},G_{2},C))=((F,G,A)\cup _{\mathcal{R}}(F_{1},G_{1},B))\cap _{\mathcal{%
R}}((F,G,A)\cup _{\mathcal{R}}(F_{2},G_{2},C))$,

\item $(F,G,A)\cup _{\mathcal{R}}((F_{1},G_{1},B)\tilde{\cap}%
(F_{2},G_{2},C))=((F,G,A)\cup _{\mathcal{R}}(F_{1},G_{1},B))\tilde{\cap}%
((F,G,A)\cup _{\mathcal{R}}(F_{2},G_{2},C))$,

\item $(F,G,A)\tilde{\cup}((F_{1},G_{1},B)\tilde{\cap}(F_{2},G_{2},C))\tilde{%
\supset}((F,G,A)\tilde{\cup}(F_{1},G_{1},B))\tilde{\cap}((F,G,A)\tilde{\cup}$
$(F_{2},G_{2},C))$,

\item $(F,G,A)\tilde{\cup}((F_{1},G_{1},B)\cup _{\mathcal{R}}(F_{2},G_{2},C))%
\tilde{\subset}((F,G,A)\tilde{\cup}(F_{1},G_{1},B))\cup _{\mathcal{R}%
}((F,G,A)\tilde{\cup}$ $(F_{2},G_{2},C))$,

\item $(F,G,A)\tilde{\cup}((F_{1},G_{1},B)\cap _{\mathcal{R}%
}(F_{2},G_{2},C))=((F,G,A)\tilde{\cup}(F_{1},G_{1},B))\cap _{\mathcal{R}%
}((F,G,A)\tilde{\cup}$ $(F_{2},G_{2},C))$,

\item $(F,G,A)\tilde{\cap}((F_{1},G_{1},B)\tilde{\cup}(F_{2},G_{2},C))\tilde{%
\subset}((F,G,A)\tilde{\cap}(F_{1},G_{1},B))\tilde{\cup}((F,G,A)\tilde{\cap}$
$(F_{2},G_{2},C))$,

\item $(F,G,A)\tilde{\cap}((F_{1},G_{1},B)\cup _{\mathcal{R}%
}(F_{2},G_{2},C))=((F,G,A)\tilde{\cap}(F_{1},G_{1},B))\cup _{\mathcal{R}%
}((F,G,A)\tilde{\cap}$ $(F_{2},G_{2},C))$,

\item $(F,G,A)\tilde{\cap}((F_{1},G_{1},B)\cap _{\mathcal{R}}(F_{2},G_{2},C))%
\tilde{\supset}((F,G,A)\tilde{\cap}(F_{1},G_{1},B))\cap _{\mathcal{R}%
}((F,G,A)\tilde{\cap}$ $(F_{2},G_{2},C))$.
\end{enumerate}

\begin{proof}
\ \ \ 

\begin{enumerate}
\item For any $e\in A\cap (B\cup C)$, we have following three disjoint cases:

$\mathbf{(i)}$ If $e\in A\cap (B\setminus C)$, then 
\begin{eqnarray*}
(F\cap _{\mathcal{R}}(F_{1}\tilde{\cup}F_{2}))(e) &=&F(e)\wedge F_{1}(e) \\
(G\cap _{\mathcal{R}}(G_{1}\tilde{\cup}G_{2}))(\lnot e) &=&G(\lnot e)\vee
G_{1}(\lnot e)
\end{eqnarray*}%
and 
\begin{eqnarray*}
((F\cap _{\mathcal{R}}F_{1})\tilde{\cup}(F\cap _{\mathcal{R}}F_{2}))(e)
&=&(F\cap _{\mathcal{R}}F_{1})(e)\vee \emptyset \\
&=&F(e)\wedge F_{1}(e) \\
((G\cap _{\mathcal{R}}G_{1})\tilde{\cup}(G\cap _{\mathcal{R}}G_{2}))(\lnot
e) &=&(G\cap _{\mathcal{R}}G_{1})(\lnot e)\wedge \mathcal{U} \\
&=&G(\lnot e)\vee G_{1}(\lnot e)\text{.}
\end{eqnarray*}%
$\mathbf{(ii)}$ If $e\in A\cap (C\setminus B)$, then 
\begin{eqnarray*}
(F\cap _{\mathcal{R}}(F_{1}\tilde{\cup}F_{2}))(e) &=&F(e)\wedge F_{2}(e) \\
(G\cap _{\mathcal{R}}(G_{1}\tilde{\cup}G_{2}))(\lnot e) &=&G(\lnot e)\vee
G_{2}(\lnot e)
\end{eqnarray*}%
and 
\begin{eqnarray*}
((F\cap _{\mathcal{R}}F_{1})\tilde{\cup}(F\cap _{\mathcal{R}}F_{2}))(e)
&=&\emptyset \vee (F\cap _{\mathcal{R}}F_{2})(e) \\
&=&F(e)\wedge F_{2}(e) \\
((G\cap _{\mathcal{R}}G_{1})\tilde{\cup}(G\cap _{\mathcal{R}}G_{2}))(\lnot
e) &=&\mathcal{U}\wedge (G\cap _{\mathcal{R}}G_{2})(\lnot e) \\
&=&G(\lnot e)\vee G_{2}(\lnot e)\text{.}
\end{eqnarray*}%
$\mathbf{(iii)}$ If $e\in A\cap (B\cap C)$, then 
\begin{eqnarray*}
(F\cap _{\mathcal{R}}(F_{1}\tilde{\cup}F_{2}))(e) &=&F(e)\wedge
(F_{1}(e)\vee F_{2}(e)) \\
(G\cap _{\mathcal{R}}(G_{1}\tilde{\cup}G_{2}))(\lnot e) &=&G(\lnot e)\vee
(G_{1}(\lnot e)\wedge G_{2}(\lnot e))
\end{eqnarray*}%
and 
\begin{eqnarray*}
((F\cap _{\mathcal{R}}F_{1})\tilde{\cup}(F\cap _{\mathcal{R}}F_{2}))(e)
&=&(F\cap _{\mathcal{R}}F_{1})(e)\vee (F\cap _{\mathcal{R}}F_{2})(e) \\
&=&(F(e)\wedge F_{1}(e))\vee (F(e)\wedge F_{2}(e)) \\
&=&F(e)\wedge (F_{1}(e)\vee F_{2}(e)) \\
((G\cap _{\mathcal{R}}G_{1})\tilde{\cup}(G\cap _{\mathcal{R}}G_{2}))(\lnot
e) &=&(G\cap _{\mathcal{R}}G_{1})(\lnot e)\wedge (G\cap _{\mathcal{R}%
}G_{2})(\lnot e) \\
&=&(G(\lnot e)\vee G_{1}(\lnot e))\wedge (G(\lnot e)\vee G_{2}(\lnot e)) \\
&=&G(\lnot e)\vee (G_{1}(\lnot e)\wedge G_{2}(\lnot e))\text{.}
\end{eqnarray*}%
Thus

$(F,G,A)\cap _{\mathcal{R}}((F_{1},G_{1},B)\tilde{\cup}%
(F_{2},G_{2},C))=((F,G,A)\cap _{\mathcal{R}}(F_{1},G_{1},B))\tilde{\cup}%
((F,G,A)\cap _{\mathcal{R}}(F_{2},G_{2},C))$.
\end{enumerate}

Similarly, we can check for the remaining parts.
\end{proof}
\end{proposition}

Now we consider the collection of all bipolar soft sets over $U$ and denote
it by $\mathcal{BSS}(U)^{E}$ and let us denote its sub collection of all
bipolar soft sets over $U$ with fixed set of parameters $A$ by $\mathcal{BSS}%
(U)_{A}$. We note that this collection is partially ordered by inclusion. We
conclude from above results that:

\begin{proposition}
$(\mathcal{BSS}(U)^{E},\tilde{\cap},\cup _{\mathcal{R}})$ and $(\mathcal{BSS}%
(U)^{E},\tilde{\cup},\cap _{\mathcal{R}})$ are distributive lattices and $(%
\mathcal{BSS}(U)^{E},\cup _{\mathcal{R}},\tilde{\cap})$ and $(\mathcal{BSS}%
(U)^{E},\cap _{\mathcal{R}},\tilde{\cup})$ are their duals respectively.

\begin{proof}
Follows from above results.
\end{proof}
\end{proposition}

\begin{proposition}
$(\mathcal{BSS}(U)^{E},\cap _{\mathcal{R}},\tilde{\cup})$\ is a bounded
distributive lattice, with least element $(\Phi ,\mathfrak{U},\varnothing )$
and greatest element $(\mathfrak{U},\Phi ,E),$ while $(\mathcal{BSS}(U)^{E},%
\tilde{\cup},\cap _{\mathcal{R}},(\mathfrak{U},\Phi ,E),(\Phi ,\mathfrak{U}%
,\varnothing ))$ is its dual.

\begin{proof}
Follows from above results.
\end{proof}
\end{proposition}

\begin{proposition}
$(\mathcal{BSS}(U)_{A},\cap _{\mathcal{R}},\tilde{\cup})=(\mathcal{BSS}%
(U)_{A},\tilde{\cap},\cup _{\mathcal{R}})$ is a bounded distributive
lattice, with least element $(\Phi ,\mathfrak{U},A)$ and greatest element $(%
\mathfrak{U},\Phi ,A)$.

\begin{proof}
Follows from above results.
\end{proof}
\end{proposition}

\begin{proposition}
Let $(F,$ $G,$ $A)$ and $(F_{1},G_{1},A)$ be two \textit{bipolar soft set}s
over a common universe $U$. Then

\begin{enumerate}
\item $((F,G,A)^{c})^{c}=(F,$ $G,$ $A)$,

\item $(F,G,A)\tilde{\subseteq}(F_{1},G_{1},A)$ implies $(F_{1},G_{1},A)^{c}%
\tilde{\subseteq}(F,G,A)^{c}$.
\end{enumerate}

\begin{proof}
1. is straightforward

2. If $(F,G,A)\tilde{\subseteq}(F_{1},G_{1},A)$ then%
\begin{equation*}
F(e)\subseteq F_{1}(e)\text{ and }G_{1}(\lnot e)\subseteq G(\lnot e)\text{
for all }e\in A
\end{equation*}%
implies that%
\begin{equation*}
(G_{1},F_{1},A)\tilde{\subseteq}(G,F,A)\text{.}
\end{equation*}

Hence $(F_{1},G_{1},A)^{c}\tilde{\subseteq}(F,G,A)^{c}$.
\end{proof}
\end{proposition}

\begin{proposition}
$(\mathcal{BSS}(U)_{A},\cap _{\mathcal{R}},\cup _{\mathcal{R}},^{c},(%
\mathfrak{U},\Phi ,A),(\Phi ,\mathfrak{U},A))$ is\ a De Morgan algebra.

\begin{proof}
Straightforward.
\end{proof}
\end{proposition}

\section{Application of Bipolar Soft Sets in a Decision Making Problem}

Decision making is an important factor of all scientific professions where
experts apply their knowledge in that area to make decisions wisely. We
apply the concept of bipolar soft sets for modelling of a given problem and
then we give an algorithm for the choice of optimal object based upon the
available sets of information. Let $U$ be the initial universe and $E$\ be a
set of parameters.

For the data analysis of a bipolar soft set, we shall use the single table
representation of $(F,G,E)$ as discussed on Page \pageref{tbldesc1}. It is
understood that both types of tabular presentations are equivalent and may
be used interchangeably where required. We shall adapt the following
terminology afterwards:

\begin{definition}
The \textit{decision value} of an object $m_{i}\in U$ is $d_{i}$, given by%
\begin{equation*}
d_{i}=\underset{j}{\dsum }a_{ij}
\end{equation*}%
where $a_{ij}$ is $(i,j)-th$ entry in the table of the bipolar soft set. We
adjoin the column of decision parameter $d$ having values $d_{i}$, with the
table of bipolar soft set $(F,G,E)$ to obtain the \textit{decision table}.
\end{definition}

We define the concept of indiscernibility relations associated with a
bipolar soft set.

\begin{definition}
If $(F,G,E)$ is a bipolar soft set over $U$ along with the set $E$ of choice
parameters, then:
\end{definition}

\begin{enumerate}
\item $\emptyset \neq F(e)\subset U$, $\emptyset \neq G(\lnot e)\subset U$,
such that $F(e)\cup G(\lnot e)\neq U$, partition $U$ into three classes;

\item If any one of $F(e)$ and $G(\lnot e)$ is empty and the other one is a
proper subset of $U$, then we have two classes of elements in $U$;

\item If any one of $F(e)$ and $G(\lnot e)$ is equal to $U$, it provides the
universal equivalence relation $U\times U$.
\end{enumerate}

In either of above three cases, these classes correspond to an equivalence
relation on $U$. Consequently, we see that for each parameter $e\in E$, we
have an equivalence relation on $U$. If this equivalence relation is denoted
by $\sigma (e)$ for all $e\in E$, then $(\sigma ,E)$ is a soft equivalence
relation over $U$. We denote 
\begin{equation*}
IND(F,G,E)=\underset{e\in E}{\cap }\sigma (e)\text{.}
\end{equation*}%
Clearly $IND(F,G,E)$ is itself an equivalence relation on $U$. The classes
of $IND(F,G,E)$ are basic categories of knowledge presented by a bipolar
soft set over $U$. We may further consider, that, $IND(E)=IND(F,G,E)$ where $%
E$ is the set of parameters. We shall say that a decision table of $(F,G,E)$
is consistent if and only if $IND(E)\subseteq IND(D)$, where $IND(D)$ is the
equivalence relation that classifies $U$ into the categories having the same
decision values.

\begin{definition}
\label{defcore}Let $T$ be a consistent decision table of bipolar soft set $%
(F,G,E)$ and $T_{\gamma }$ be a decision table obtained from $T$ by
eliminating some column of $\gamma \in C$. Then $\gamma $ is dispensable in $%
T$ if

\begin{enumerate}
\item $T_{\gamma }$ is consistent that is $IND(C-\gamma )=IND(D_{\gamma })$

\item $IND(D)=IND(D_{\gamma })$
\end{enumerate}
\end{definition}

Otherwise $\gamma $\ is indispensable or core parameter. The set of all core
parameters of $C$ is denoted by $CORE(C)$.

\begin{algorithm}
The algorithm for the selection of the best choice is given as:

\begin{enumerate}
\item Input the bipolar soft set $(F,G,E)$.

\item Input the set of choice parameters $C\subseteq E$.

\item Input the decision parameter $d\in D$, $d_{i}=\underset{j}{\dsum }%
a_{ij}$ as the last column in the table obtained by choice parameters.

\item Rearrange the Input by placing the objects having the same value for
the parameter $d$ adjacent to each other.

\item Distinguish the objects with different values of $d$ by double line.

\item Identify core parameters as defined in Definition \ref{defcore}.
Eliminate all the dispensable parameters one by one, resulting a table with
minimum number of condition parameters having the same classification
ability for $d$ as the original table with $d$.

\item Find $k$, for which $d_{k}=\max d_{i}$.
\end{enumerate}

Then $m_{k}$ is the optimal choice object. If $k$ has more than one values,
then any one of $m_{k}$'s can be chosen.
\end{algorithm}

\begin{example}
Let $U=\{m_{1},m_{2},m_{3},m_{4},m_{5},m_{6},m_{7},m_{8}\}$ be the set of
candidates who have applied for a job position of Office Representative in
Customer Care Centre of a company. Let $E=%
\{e_{1},e_{2},e_{3},e_{4},e_{5},e_{6},e_{7},e_{8},e_{9}\}=\{$Hard Working,
Optimism, Enthusiasm, Individualism, Imaginative, Flexibility, Decisiveness,
Self-confidence, Politeness$\}$ and $\lnot E=\{\lnot e_{1},\lnot e_{2},\lnot
e_{3},\lnot e_{4},\lnot e_{5},\lnot e_{6},\lnot e_{7},\lnot e_{8},\lnot
e_{9}\}=\{$Negligent, Pessimism, Half-hearted, Dependence, Unimaginative,
Rigidity, Indecisiveness, Shyness, Harshness$\ \}$. Here the gray area is
obviously the moderate form of parameters. Let the bipolar soft sets $(F$,$G$%
,$E)$ describes the \textquotedblleft\ Personality Analysis of
Candidates\textquotedblright\ as:%
\begin{eqnarray*}
F(e_{1})
&=&\{m_{1},m_{4},m_{5},m_{8}\},F(e_{2})=\{m_{1},m_{2},m_{3},m_{4},m_{8}\}, \\
F(e_{3})
&=&\{m_{2},m_{4},m_{6},m_{7},m_{8}\},F(e_{4})=\{m_{6},m_{7}\},F(e_{5})=%
\{m_{1},m_{7},m_{8}\}, \\
F(e_{6})
&=&\{m_{4},m_{5},m_{6},m_{7}\},F(e_{7})=\{m_{1},m_{2},m_{5},m_{6},m_{8}%
\},F(e_{8})=\{m_{1},m_{6},m_{8}\}, \\
F(e_{9}) &=&\{m_{2},m_{3},m_{4},m_{6},m_{7}\}\text{.}
\end{eqnarray*}%
and%
\begin{eqnarray*}
G(\lnot e_{1}) &=&\{m_{6},m_{7}\},G(\lnot e_{2})=\{m_{5},m_{6}\},G(\lnot
e_{3})=\{\},G(\lnot e_{4})=\{m_{1},m_{3},m_{8}\}, \\
G(\lnot e_{5}) &=&\{m_{2},m_{3},m_{4},m_{5},m_{6}\},G(\lnot
e_{6})=\{m_{8}\},G(\lnot e_{7})=\{m_{3},m_{4}\}, \\
G(\lnot e_{8}) &=&\{m_{5}\},G(\lnot e_{9})=\{m_{1},m_{5}\}\text{.}
\end{eqnarray*}

\begin{enumerate}
\item Input the bipolar soft set $(F,G,E)$ given by Table \ref{TableKey
copy(3)}.

\begin{table}[h] \centering%
\begin{tabular}{|l|l|l|l|l|l|l|l|l|l|}
\hline
$(F,G,E)$ & $e_{1}$ & $e_{2}$ & $e_{3}$ & $e_{4}$ & $e_{5}$ & $e_{6}$ & $%
e_{7}$ & $e_{8}$ & $e_{9}$ \\ \hline
$m_{1}$ & \multicolumn{1}{|r|}{$1$} & \multicolumn{1}{|r|}{$1$} & 
\multicolumn{1}{|r|}{$0$} & \multicolumn{1}{|r|}{$-1$} & 
\multicolumn{1}{|r|}{$1$} & \multicolumn{1}{|r|}{$0$} & \multicolumn{1}{|r|}{%
$1$} & \multicolumn{1}{|r|}{$1$} & \multicolumn{1}{|r|}{$-1$} \\ \hline
$m_{2}$ & \multicolumn{1}{|r|}{$0$} & \multicolumn{1}{|r|}{$1$} & 
\multicolumn{1}{|r|}{$1$} & \multicolumn{1}{|r|}{$0$} & \multicolumn{1}{|r|}{%
$-1$} & \multicolumn{1}{|r|}{$0$} & \multicolumn{1}{|r|}{$1$} & 
\multicolumn{1}{|r|}{$0$} & \multicolumn{1}{|r|}{$1$} \\ \hline
$m_{3}$ & \multicolumn{1}{|r|}{$0$} & \multicolumn{1}{|r|}{$1$} & 
\multicolumn{1}{|r|}{$0$} & \multicolumn{1}{|r|}{$-1$} & 
\multicolumn{1}{|r|}{$-1$} & \multicolumn{1}{|r|}{$0$} & 
\multicolumn{1}{|r|}{$-1$} & \multicolumn{1}{|r|}{$0$} & 
\multicolumn{1}{|r|}{$1$} \\ \hline
$m_{4}$ & \multicolumn{1}{|r|}{$1$} & \multicolumn{1}{|r|}{$1$} & 
\multicolumn{1}{|r|}{$1$} & \multicolumn{1}{|r|}{$0$} & \multicolumn{1}{|r|}{%
$-1$} & \multicolumn{1}{|r|}{$1$} & \multicolumn{1}{|r|}{$-1$} & 
\multicolumn{1}{|r|}{$0$} & \multicolumn{1}{|r|}{$1$} \\ \hline
$m_{5}$ & \multicolumn{1}{|r|}{$1$} & \multicolumn{1}{|r|}{$-1$} & 
\multicolumn{1}{|r|}{$0$} & \multicolumn{1}{|r|}{$0$} & \multicolumn{1}{|r|}{%
$-1$} & \multicolumn{1}{|r|}{$1$} & \multicolumn{1}{|r|}{$1$} & 
\multicolumn{1}{|r|}{$-1$} & \multicolumn{1}{|r|}{$-1$} \\ \hline
$m_{6}$ & \multicolumn{1}{|r|}{$-1$} & \multicolumn{1}{|r|}{$-1$} & 
\multicolumn{1}{|r|}{$1$} & \multicolumn{1}{|r|}{$1$} & \multicolumn{1}{|r|}{%
$-1$} & \multicolumn{1}{|r|}{$1$} & \multicolumn{1}{|r|}{$1$} & 
\multicolumn{1}{|r|}{$1$} & \multicolumn{1}{|r|}{$1$} \\ \hline
$m_{7}$ & \multicolumn{1}{|r|}{$-1$} & \multicolumn{1}{|r|}{$0$} & 
\multicolumn{1}{|r|}{$1$} & \multicolumn{1}{|r|}{$1$} & \multicolumn{1}{|r|}{%
$1$} & \multicolumn{1}{|r|}{$1$} & \multicolumn{1}{|r|}{$0$} & 
\multicolumn{1}{|r|}{$0$} & \multicolumn{1}{|r|}{$1$} \\ \hline
$m_{8}$ & \multicolumn{1}{|r|}{$1$} & \multicolumn{1}{|r|}{$1$} & 
\multicolumn{1}{|r|}{$1$} & \multicolumn{1}{|r|}{$-1$} & 
\multicolumn{1}{|r|}{$1$} & \multicolumn{1}{|r|}{$-1$} & 
\multicolumn{1}{|r|}{$1$} & \multicolumn{1}{|r|}{$1$} & \multicolumn{1}{|r|}{%
$0$} \\ \hline
\end{tabular}%
\caption{Table of (F.G,E)}\label{TableKey copy(3)}%
\end{table}%

\item Let $C=\{e_{1},e_{3},e_{4},e_{5},e_{7},e_{8}\}$.

\item The decision table of bipolar soft set $(F,G,C)$ is given by Table \ref%
{TableKey copy(4)}.

\begin{table}[h] \centering%
\begin{tabular}{|l|l|l|l|l|l|l|l|}
\hline
$(F,G,C)$ & $e_{1}$ & $e_{3}$ & $e_{4}$ & $e_{5}$ & $e_{7}$ & $e_{8}$ & $d$
\\ \hline
$m_{1}$ & \multicolumn{1}{|r|}{$1$} & \multicolumn{1}{|r|}{$0$} & 
\multicolumn{1}{|r|}{$-1$} & \multicolumn{1}{|r|}{$1$} & 
\multicolumn{1}{|r|}{$1$} & \multicolumn{1}{|r|}{$1$} & \multicolumn{1}{|r|}{%
$3$} \\ \hline
$m_{2}$ & \multicolumn{1}{|r|}{$0$} & \multicolumn{1}{|r|}{$1$} & 
\multicolumn{1}{|r|}{$0$} & \multicolumn{1}{|r|}{$-1$} & 
\multicolumn{1}{|r|}{$1$} & \multicolumn{1}{|r|}{$0$} & \multicolumn{1}{|r|}{%
$1$} \\ \hline
$m_{3}$ & \multicolumn{1}{|r|}{$0$} & \multicolumn{1}{|r|}{$0$} & 
\multicolumn{1}{|r|}{$-1$} & \multicolumn{1}{|r|}{$-1$} & 
\multicolumn{1}{|r|}{$-1$} & \multicolumn{1}{|r|}{$0$} & 
\multicolumn{1}{|r|}{$-3$} \\ \hline
$m_{4}$ & \multicolumn{1}{|r|}{$1$} & \multicolumn{1}{|r|}{$1$} & 
\multicolumn{1}{|r|}{$0$} & \multicolumn{1}{|r|}{$-1$} & 
\multicolumn{1}{|r|}{$-1$} & \multicolumn{1}{|r|}{$0$} & 
\multicolumn{1}{|r|}{$0$} \\ \hline
$m_{5}$ & \multicolumn{1}{|r|}{$1$} & \multicolumn{1}{|r|}{$0$} & 
\multicolumn{1}{|r|}{$0$} & \multicolumn{1}{|r|}{$-1$} & 
\multicolumn{1}{|r|}{$1$} & \multicolumn{1}{|r|}{$-1$} & 
\multicolumn{1}{|r|}{$0$} \\ \hline
$m_{6}$ & \multicolumn{1}{|r|}{$-1$} & \multicolumn{1}{|r|}{$1$} & 
\multicolumn{1}{|r|}{$1$} & \multicolumn{1}{|r|}{$-1$} & 
\multicolumn{1}{|r|}{$1$} & \multicolumn{1}{|r|}{$1$} & \multicolumn{1}{|r|}{%
$2$} \\ \hline
$m_{7}$ & \multicolumn{1}{|r|}{$-1$} & \multicolumn{1}{|r|}{$1$} & 
\multicolumn{1}{|r|}{$1$} & \multicolumn{1}{|r|}{$1$} & \multicolumn{1}{|r|}{%
$0$} & \multicolumn{1}{|r|}{$0$} & \multicolumn{1}{|r|}{$2$} \\ \hline
$m_{8}$ & \multicolumn{1}{|r|}{$1$} & \multicolumn{1}{|r|}{$1$} & 
\multicolumn{1}{|r|}{$-1$} & \multicolumn{1}{|r|}{$1$} & 
\multicolumn{1}{|r|}{$1$} & \multicolumn{1}{|r|}{$1$} & \multicolumn{1}{|r|}{%
$4$} \\ \hline
\end{tabular}%
\caption{Decision Table for (F,G,C)}\label{TableKey copy(4)}%
\end{table}%

We note that%
\begin{eqnarray*}
IND(C)
&=&%
\{(m_{1},m_{1}),(m_{2},m_{2}),(m_{3},m_{3}),(m_{4},m_{4}),(m_{5},m_{5}),(m_{6},m_{6}),(m_{7},m_{7}),(m_{8},m_{8})\}
\\
&\subset
&\{(m_{1},m_{1}),(m_{2},m_{2}),(m_{3},m_{3}),(m_{4},m_{4}),(m_{5},m_{5}),(m_{6},m_{6}),(m_{7},m_{7}),(m_{8},m_{8}),
\\
&&(m_{5},m_{4}),(m_{4},m_{5}),(m_{6},m_{7}),(m_{6},m_{7})\} \\
&=&IND(D)\text{.}
\end{eqnarray*}%
Hence the decision table is consistent.

\item We rearrange the table according to the same values for $d$ to obtain
Table \ref{TableKey copy(5)}.

\begin{table}[h] \centering%
\begin{tabular}{||l||l|l|l|l|l|l|l|}
\hline\hline
$(F,G,C)$ & $e_{1}$ & \multicolumn{1}{||l|}{$e_{3}$} & \multicolumn{1}{||l|}{%
$e_{4}$} & \multicolumn{1}{||l|}{$e_{5}$} & \multicolumn{1}{||l|}{$e_{7}$} & 
\multicolumn{1}{||l|}{$e_{8}$} & \multicolumn{1}{||l||}{$d$} \\ \hline\hline
$m_{8}$ & \multicolumn{1}{||r|}{$1$} & \multicolumn{1}{|r|}{$1$} & 
\multicolumn{1}{|r|}{$-1$} & \multicolumn{1}{|r|}{$1$} & 
\multicolumn{1}{|r|}{$1$} & \multicolumn{1}{|r|}{$1$} & \multicolumn{1}{|r|}{%
$4$} \\ \hline\hline
$m_{1}$ & \multicolumn{1}{||r|}{$1$} & \multicolumn{1}{|r|}{$0$} & 
\multicolumn{1}{|r|}{$-1$} & \multicolumn{1}{|r|}{$1$} & 
\multicolumn{1}{|r|}{$1$} & \multicolumn{1}{|r|}{$1$} & \multicolumn{1}{|r|}{%
$3$} \\ \hline\hline
$m_{6}$ & \multicolumn{1}{||r|}{$-1$} & \multicolumn{1}{|r|}{$1$} & 
\multicolumn{1}{|r|}{$1$} & \multicolumn{1}{|r|}{$-1$} & 
\multicolumn{1}{|r|}{$1$} & \multicolumn{1}{|r|}{$1$} & \multicolumn{1}{|r|}{%
$2$} \\ \hline
$m_{7}$ & \multicolumn{1}{||r|}{$-1$} & \multicolumn{1}{|r|}{$1$} & 
\multicolumn{1}{|r|}{$1$} & \multicolumn{1}{|r|}{$1$} & \multicolumn{1}{|r|}{%
$0$} & \multicolumn{1}{|r|}{$0$} & \multicolumn{1}{|r|}{$2$} \\ \hline\hline
$m_{2}$ & \multicolumn{1}{||r|}{$0$} & \multicolumn{1}{|r|}{$1$} & 
\multicolumn{1}{|r|}{$0$} & \multicolumn{1}{|r|}{$-1$} & 
\multicolumn{1}{|r|}{$1$} & \multicolumn{1}{|r|}{$0$} & \multicolumn{1}{|r|}{%
$1$} \\ \hline\hline
$m_{4}$ & \multicolumn{1}{||r|}{$1$} & \multicolumn{1}{|r|}{$1$} & 
\multicolumn{1}{|r|}{$0$} & \multicolumn{1}{|r|}{$-1$} & 
\multicolumn{1}{|r|}{$-1$} & \multicolumn{1}{|r|}{$0$} & 
\multicolumn{1}{|r|}{$0$} \\ \hline
$m_{5}$ & \multicolumn{1}{||r|}{$1$} & \multicolumn{1}{|r|}{$0$} & 
\multicolumn{1}{|r|}{$0$} & \multicolumn{1}{|r|}{$-1$} & 
\multicolumn{1}{|r|}{$1$} & \multicolumn{1}{|r|}{$-1$} & 
\multicolumn{1}{|r|}{$0$} \\ \hline\hline
$m_{3}$ & \multicolumn{1}{||r|}{$0$} & \multicolumn{1}{|r|}{$0$} & 
\multicolumn{1}{|r|}{$-1$} & \multicolumn{1}{|r|}{$-1$} & 
\multicolumn{1}{|r|}{$-1$} & \multicolumn{1}{|r|}{$0$} & 
\multicolumn{1}{|r|}{$-3$} \\ \hline
\end{tabular}%
\caption{Table obtained after rearrangement of Rows according to the values
of d}\label{TableKey copy(5)}%
\end{table}%

\item We identify that%
\begin{equation*}
CORE(C)=C\text{.}
\end{equation*}

\item We conclude from the values of $d$ that $d_{8}=\max d_{i}=4$ and hence 
$k=8$.
\end{enumerate}

Thus $m_{8}$ is the optimal choice object and so $m_{8}$ is the best
candidate for the position. In case that $m_{8}$ can not join the position, $%
m_{1}$ will be selected, and if $m_{1}$ will also not be able to join then
either $m_{6}$ or $m_{7}$ may be selected.
\end{example}

We define another type of weighted table of the reduct bipolar soft set $%
(F,G,C)$. The motivation is, that, some of the parameters are of less
importance than the other ones so they must be graded with lesser priority.
For this reason, we suggest that the column of that parameter will have
entries 
\begin{equation*}
b_{ij}=\left\{ 
\begin{array}{lll}
a_{ij}\times w_{j} &  & \text{if }a_{ij}=1 \\ 
0 &  & \text{if }a_{ij}=0 \\ 
a_{ij}\times (1-w_{j}) &  & \text{if }a_{ij}=-1%
\end{array}%
\right. ,
\end{equation*}%
instead of $0$ and $1$ and $-1$ only, where $a_{ij}$ are the entries in the
table of the reduct bipolar soft set $(F,G,C)$.

\begin{definition}
The \textit{weighted decision value} of an object $m_{i}\in U$ is 
\begin{equation*}
d_{i}=\underset{j}{\dsum }b_{ij}\text{.}
\end{equation*}
\end{definition}

The revised algorithm will be given now:

\begin{algorithm}
The algorithm for the selection of the best choice is given as:

\begin{enumerate}
\item Input the bipolar soft set $(F,G,E)$.

\item Input the set of choice parameters $C\subseteq E$.

\item Find weighted table of the bipolar soft set $(F,G,C)$ according to the
weights decided.

\item Input the decision parameter $d\in D$, $d_{i}=\underset{j}{\dsum }%
b_{ij}$ as the last column in the weighted table $T_{w}$.

\item Rearrange the input by placing the objects having the same value for
the parameter $d$ adjacent to each other.

\item Distinguish the objects with different values of $d$ by double line.

\item Identify core parameters. Eliminate all the dispensable parameters one
by one, resulting a table with minimum number of condition parameters having
the same classification ability for $d$ as the original table with $d$.

\item Find $k$, for which $d_{k}=\max d_{i}$.
\end{enumerate}
\end{algorithm}

Then $m_{k}$ is the optimal choice object. If $k$ has more than one values,
then any one of $m_{k}$'s can be chosen.

Now we solve the original problem using this revised algorithm. Suppose that
the selection board sets the following weights for parameters of $C$ and
take start from the 3rd step as:%
\begin{eqnarray*}
e_{1} &:&w_{1}=0.9 \\
e_{3} &:&w_{3}=0.8 \\
e_{4} &:&w_{4}=0.5 \\
e_{5} &:&w_{5}=0.6 \\
e_{7} &:&w_{7}=0.9 \\
e_{8} &:&w_{8}=0.9
\end{eqnarray*}%
The weighted decision table of bipolar soft set $(F,G,C)$ is given by Table %
\ref{TableKey copy(8)}.

\begin{table}[h] \centering%
\begin{tabular}{|l|l|l|l|l|l|l|l|}
\hline
$(F,G,C)_{w}$ & $e_{1}$ & $e_{3}$ & $e_{4}$ & $e_{5}$ & $e_{7}$ & $e_{8}$ & $%
d$ \\ \hline
$m_{1}$ & \multicolumn{1}{|r|}{$0.9$} & \multicolumn{1}{|r|}{$0$} & 
\multicolumn{1}{|r|}{$-0.5$} & \multicolumn{1}{|r|}{$0.6$} & 
\multicolumn{1}{|r|}{$0.9$} & \multicolumn{1}{|r|}{$0.9$} & 
\multicolumn{1}{|r|}{$2.8$} \\ \hline
$m_{2}$ & \multicolumn{1}{|r|}{$0$} & \multicolumn{1}{|r|}{$0.8$} & 
\multicolumn{1}{|r|}{$0$} & \multicolumn{1}{|r|}{$-0.4$} & 
\multicolumn{1}{|r|}{$0.9$} & \multicolumn{1}{|r|}{$0$} & 
\multicolumn{1}{|r|}{$1.3$} \\ \hline
$m_{3}$ & \multicolumn{1}{|r|}{$0$} & \multicolumn{1}{|r|}{$0$} & 
\multicolumn{1}{|r|}{$-0.5$} & \multicolumn{1}{|r|}{$-0.4$} & 
\multicolumn{1}{|r|}{$-0.1$} & \multicolumn{1}{|r|}{$0$} & 
\multicolumn{1}{|r|}{$-1$} \\ \hline
$m_{4}$ & \multicolumn{1}{|r|}{$0.9$} & \multicolumn{1}{|r|}{$0.8$} & 
\multicolumn{1}{|r|}{$0$} & \multicolumn{1}{|r|}{$-0.4$} & 
\multicolumn{1}{|r|}{$-0.1$} & \multicolumn{1}{|r|}{$0$} & 
\multicolumn{1}{|r|}{$1.2$} \\ \hline
$m_{5}$ & \multicolumn{1}{|r|}{$0.9$} & \multicolumn{1}{|r|}{$0.8$} & 
\multicolumn{1}{|r|}{$0$} & \multicolumn{1}{|r|}{$-0.4$} & 
\multicolumn{1}{|r|}{$0.9$} & \multicolumn{1}{|r|}{$-0.1$} & 
\multicolumn{1}{|r|}{$2.1$} \\ \hline
$m_{6}$ & \multicolumn{1}{|r|}{$-0.1$} & \multicolumn{1}{|r|}{$0.8$} & 
\multicolumn{1}{|r|}{$0.5$} & \multicolumn{1}{|r|}{$-0.4$} & 
\multicolumn{1}{|r|}{$0.9$} & \multicolumn{1}{|r|}{$0.9$} & 
\multicolumn{1}{|r|}{$2.6$} \\ \hline
$m_{7}$ & \multicolumn{1}{|r|}{$-0.1$} & \multicolumn{1}{|r|}{$0.8$} & 
\multicolumn{1}{|r|}{$0.5$} & \multicolumn{1}{|r|}{$0.6$} & 
\multicolumn{1}{|r|}{$0$} & \multicolumn{1}{|r|}{$0$} & \multicolumn{1}{|r|}{%
$1.8$} \\ \hline
$m_{8}$ & \multicolumn{1}{|r|}{$0.9$} & \multicolumn{1}{|r|}{$0.8$} & 
\multicolumn{1}{|r|}{$-0.5$} & \multicolumn{1}{|r|}{$0.6$} & 
\multicolumn{1}{|r|}{$0.9$} & \multicolumn{1}{|r|}{$0.9$} & 
\multicolumn{1}{|r|}{$3.6$} \\ \hline
\end{tabular}%
\caption{Weighted Decision Table for (F,G,C)}\label{TableKey copy(8)}%
\end{table}%

We note that%
\begin{eqnarray*}
IND(C)
&=&%
\{(m_{1},m_{1}),(m_{2},m_{2}),(m_{3},m_{3}),(m_{4},m_{4}),(m_{5},m_{5}),(m_{6},m_{6}),(m_{7},m_{7}),(m_{8},m_{8})\}
\\
&=&IND(D)\text{.}
\end{eqnarray*}%
Hence the decision table is consistent. We rearrange the table according to
the descending values for $d$ to obtain Table \ref{TableKey copy(6)}.

\begin{table}[h] \centering%
\begin{tabular}{||l||l|l|l|l|l|l|l|}
\hline\hline
$(F,G,C)_{w}$ & $e_{1}$ & \multicolumn{1}{||l|}{$e_{3}$} & 
\multicolumn{1}{||l|}{$e_{4}$} & \multicolumn{1}{||l|}{$e_{5}$} & 
\multicolumn{1}{||l|}{$e_{7}$} & \multicolumn{1}{||l|}{$e_{8}$} & 
\multicolumn{1}{||l||}{$d$} \\ \hline\hline
$m_{8}$ & \multicolumn{1}{||r|}{$0.9$} & \multicolumn{1}{|r|}{$0.8$} & 
\multicolumn{1}{|r|}{$-0.5$} & \multicolumn{1}{|r|}{$0.6$} & 
\multicolumn{1}{|r|}{$0.9$} & \multicolumn{1}{|r|}{$0.9$} & 
\multicolumn{1}{|r|}{$3.6$} \\ \hline\hline
$m_{1}$ & \multicolumn{1}{||r|}{$0.9$} & \multicolumn{1}{|r|}{$0$} & 
\multicolumn{1}{|r|}{$-0.5$} & \multicolumn{1}{|r|}{$0.6$} & 
\multicolumn{1}{|r|}{$0.9$} & \multicolumn{1}{|r|}{$0.9$} & 
\multicolumn{1}{|r|}{$2.8$} \\ \hline\hline
$m_{6}$ & \multicolumn{1}{||r|}{$-0.1$} & \multicolumn{1}{|r|}{$0.8$} & 
\multicolumn{1}{|r|}{$0.5$} & \multicolumn{1}{|r|}{$-0.4$} & 
\multicolumn{1}{|r|}{$0.9$} & \multicolumn{1}{|r|}{$0.9$} & 
\multicolumn{1}{|r|}{$2.6$} \\ \hline\hline
$m_{5}$ & \multicolumn{1}{||r|}{$0.9$} & \multicolumn{1}{|r|}{$0.8$} & 
\multicolumn{1}{|r|}{$0$} & \multicolumn{1}{|r|}{$-0.4$} & 
\multicolumn{1}{|r|}{$0.9$} & \multicolumn{1}{|r|}{$-0.1$} & 
\multicolumn{1}{|r|}{$2.1$} \\ \hline\hline
$m_{7}$ & \multicolumn{1}{||r|}{$-0.1$} & \multicolumn{1}{|r|}{$0.8$} & 
\multicolumn{1}{|r|}{$0.5$} & \multicolumn{1}{|r|}{$0.6$} & 
\multicolumn{1}{|r|}{$0$} & \multicolumn{1}{|r|}{$0$} & \multicolumn{1}{|r|}{%
$1.8$} \\ \hline\hline
$m_{2}$ & \multicolumn{1}{||r|}{$0$} & \multicolumn{1}{|r|}{$0.8$} & 
\multicolumn{1}{|r|}{$0$} & \multicolumn{1}{|r|}{$-0.4$} & 
\multicolumn{1}{|r|}{$0.9$} & \multicolumn{1}{|r|}{$0$} & 
\multicolumn{1}{|r|}{$1.3$} \\ \hline\hline
$m_{4}$ & \multicolumn{1}{||r|}{$0.9$} & \multicolumn{1}{|r|}{$0.8$} & 
\multicolumn{1}{|r|}{$0$} & \multicolumn{1}{|r|}{$-0.4$} & 
\multicolumn{1}{|r|}{$-0.1$} & \multicolumn{1}{|r|}{$0$} & 
\multicolumn{1}{|r|}{$1.2$} \\ \hline\hline
$m_{3}$ & \multicolumn{1}{||r|}{$0$} & \multicolumn{1}{|r|}{$0$} & 
\multicolumn{1}{|r|}{$-0.5$} & \multicolumn{1}{|r|}{$-0.4$} & 
\multicolumn{1}{|r|}{$-0.1$} & \multicolumn{1}{|r|}{$0$} & 
\multicolumn{1}{|r|}{$-1$} \\ \hline\hline
\end{tabular}%
\caption{Table of weighted Bipolar soft set (F,G,C) after Re-arrangement }%
\label{TableKey copy(6)}%
\end{table}%

We identify that%
\begin{equation*}
CORE(C)=C\text{.}
\end{equation*}%
The values of $d$ show that $d_{8}=\max d_{i}=4$ and hence $k=8$.

Once again, $m_{8}$ is the optimal choice object and so $m_{8}$ is the best
candidate for the position. We note that the difference occurs in the
position of $m_{5}$. In the first case $m_{5}$ was at $7$th position out of $%
8$ but under the weighted criteria $m_{5}$ takes $4$th position over all.

\begin{conclusion}
We have defined bipolar soft sets and various operations of union and
intersection for them. We have also shown that the concept of bipolar soft
sets is different and can not be subsumed by combining soft sets only. As we
have discussed the concept of bipolarity for the set of parameters used for
approximations of initial universe $U$, the idea may be extended for a
further study of tri-polar and hence multipolar soft sets. Therefore, this
paper gives an idea for the beginning of a new study for approximations of
data with uncertainties.
\end{conclusion}

\end{document}